\documentclass{amsart}
\usepackage{amssymb}
\usepackage{amsfonts}
\usepackage{amsmath}
\usepackage{lastpage}
\usepackage[legalpaper,bookmarks=true,colorlinks=true,linkcolor=blue,citecolor=blue]{hyperref}
\usepackage{graphicx}%
\usepackage{fancyhdr}
\usepackage{color}
\usepackage[mathlines]{lineno}
\usepackage{lscape}
\usepackage{epsfig}
\usepackage{natbib}
\usepackage{geometry}
\usepackage{tgbonum}
\fontfamily{qcr}\selectfont
\usepackage[]{listings}

\newtheorem{theorem}{Theorem}
\theoremstyle{plain}

\newtheorem{corollary}{Corollary}

\newtheorem{lemma}{Lemma}

\newtheorem{proposition}{Proposition}

\numberwithin{equation}{section}

\newcommand{\Bin}{\bigskip \noindent}

\newcommand{\Ni}{\noindent}

\begin{document}
\Large
\title[]{Statistical tests for the Pseudo-Lindley distribution and applications}

\author{Gane Samb Lo$^{(1)}$}
\author{Tchilabalo Abozou Kpanzou$^{(2)}$}
\author{Cheikh Mohamed Haidara$^{(3)}$}

\begin{abstract} The pseudo-Lindley distribution was introduced as a useful generalization of the Lindley distribution in Zeghdoudi and Nedjar (2016) who showed 
interesting properties of their new laws and efficiencies in modeling data in Reliability and Survival Analysis. In this paper we study the estimators of the 
pair of parameters and determine their asymptotic law from which a chi-square law is derived. From both asymptotic laws, statistical tests are built. Simulation studies on the tests conclude to their efficiency for data sizes generally used in Reliability. R codes related to statistical analysis on that  law are given in an appropriate archive repository code paper in Arxiv.\\

\noindent $^{(1)}$ Gane Samb Lo.\\
LERSTAD, Gaston Berger University, Saint-Louis, S\'en\'egal (main affiliation).\newline
LSTA, Pierre and Marie Curie University, Paris VI, France.\newline
AUST - African University of Sciences and Technology, Abuja, Nigeria\\
gane-samb.lo@edu.ugb.sn, gslo@aust.edu.ng, ganesamblo@ganesamblo.net\\
Permanent address : 1178 Evanston Dr NW T3P 0J9,Calgary, Alberta, Canada.

\bigskip \noindent $^{(2)}$ Kpanzou Tchilabalo Abozou (t.kpanzou@univkara.net).\\
University of Kara, Kara, Togo.

\bigskip \noindent $^{(3)}$ Cheikh Mohamed Haidara (chheikhh@yahoo.com).\\
LERSTAD, Gaston Berger University, Saint-Louis, Senegal.\\
%\noindent $^{\dag\dag}$ Mohammad Ahsanullah\\
%Department of Management Sciences. Rider University. Lawrenceville, New Jersey, USA\\
%Email : ahsan@rider.edu\\

%\noindent \noindent $^{\dag\dag\dag}$ Modou Ngom.\\
%LERSTAD, Gaston Berger University, Saint-Louis, S\'en\'egal (main affiliation).\newline
%modou.ngo.edu@gmail.com\\

\end{abstract}

\maketitle

\section{Introduction}\label{sec1}

\noindent The following probability distribution function (\textit{pdf}), named as the Pseudo-Lindley \textit{pdf}, 

\begin{equation} 
f(x)=f(x,\theta,\beta)=\frac{\theta(\beta-1+\theta x) e^{-\theta x}}{\beta} 1_{(x\geq 0)}, \label{glind}
\end{equation}

\Bin with parameters $\theta>0$ and $\beta>1$, has been introduced by \cite{Zeghoudi3}  as a generalization of the Lindley \textit{pdf} :

\begin{equation} 
\ell(x)=\frac{\theta^2(1+x) e^{-\theta x}}{1+\theta} 1_{(x\geq 0)}, \label{olind}
\end{equation}

\Bin \noindent in the sense that for $\beta=1+\theta$, $f(\circ)$ is identical to $\ell(\circ)$.\\

\Ni  Actually, $f$ derives from $\ell$ by a mixture of a Lindley distributed random variable and an independent $\Gamma(2,\theta)$ random variable with mixture coefficients $r_1=(\beta-1)/\beta$ and $r_2=1/\beta$, where $0<r_1, \ r_2<1$ and $r_1+r_2=1$.\\

\Ni The Lindley original distribution is an important law that has been used and still being used in Reliability, in Survival analysis and other important disciplines. Because of its original remarkable qualities, it kicked a considerable number of generalizations as pointed out by \cite{Zeghoudi3}. The current generalization
\eqref{glind} has been tested on real data and simulated data, and shows real interest in survival analysis, on the Guinean Ebola data for example (\cite{Zeghoudi3}). That interests that distribution demonstrated in real data modeling motivates us to give some asymptotic theories on it, in view of statistical tests.\\

\noindent This paper is a contribution of the statistical analysis of this new Pseudo-Lindley law, essentially on an asymptotic view. More precisely, we aim at proceeding to the complete asymptotic theory of the moments estimators by finding the joint Gaussian asymptotic law of the pair of estimators of $\alpha$ and $\beta$, to base statistical tests on the results, to validate them as accurate using simulation studies and to recommend sample sizes for the efficient use of our statistical tests.\\ 

\noindent Because of the considerable number of new laws, we want this paper to be an example of the use of the functional empirical process for deriving asymptotic laws as described in \cite{logfep-tool} and \cite{ips-wcia-ang} (Chapter 5, Subsection 4).\\

\noindent The remainder of the paper is organized as follows: In Section \ref{sec_02}, we consider the main tool and context. We give a review on the functional empirical process and the related asymptotic laws in Section \ref{ssec_02_01}, then the moment estimators of the parameters in Section \ref{ssec_02_01}. The joint Gaussian asymptotic law of these estimators is given in Section \ref{ssec_03}. In Section \ref{ssec_04}, a simulation study is undertaken and showed that the estimators have good performance with respect to the mean value errors, and are excellent for statistical tests even for small sample sizes. The paper is concluded in Section \ref{ssec_05}.

\section{Main tool and context} \label{sec_02}

\Ni Let us begin by the main asymptotic tool we use. \\

\subsection{Functional empirical process} \label{ssec_02_01}

\noindent  Suppose that $X$, $X_1$, $X_2$, $\cdots$ are independent and identically real-valued random variables, with common cumulative distribution function (\textit{cdf}) $F$ defined on the same probability space $(\Omega, \mathcal{A}, \mathbb{P})$. Define for each $n\geq 1,$ the
functional empirical process by 
\begin{equation*}
\mathbb{G}_{n}(f)=\frac{1}{\sqrt{n}}\sum_{i=1}^{n}(f(X_{i})-\mathbb{E}%
f(Z_{i})),
\end{equation*}

\bigskip \noindent where $f$ is a real and measurable function defined on $%
\mathbb{R}$ such that

\begin{equation}
\mathbb{V}_{X}(f)=\int \left( f(x)-\mathbb{P}_{X}(f)\right)
^{2}dP_{Z}(x)<\infty ,  \label{var}
\end{equation}

\Bin which entails

\begin{equation}
\mathbb{P}_{Z}(\left\vert f\right\vert )=\int \left\vert f(x)\right\vert
dP_{Z}(x)<\infty \text{.}  \label{esp}
\end{equation}

\bigskip \noindent Denote by $\mathcal{F}(S)$ - $\mathcal{F}$ for short -
the class of real-valued measurable functions that are defined on S such
that (\ref{var}) holds. The space $\mathcal{F}$ , when endowed with the
addition and the external multiplication by real scalars, is a linear space.
Next, it is remarkable that $\mathbb{G}_{n}$ is linear on $\mathcal{F}$, that
is for $f$ and $g$ in $\mathcal{F}$ and for $(a,b)\in \mathbb{R}{^{2}}$, we
have

\begin{equation*}
a\mathbb{G}_{n}(f)+b\mathbb{G}_{n}(g)=\mathbb{G}_{n}(af+bg).
\end{equation*}

\bigskip \noindent We have the following result.

\begin{lemma} \label{lemma2}
\bigskip Given the notation above, for any finite number of elements $%
f_{1},...,f_{k}$ of $\mathcal{S},k\geq 1,$ we have

\begin{equation*}
^{t}(\mathbb{G}_{n}(f_{1}),...,\mathbb{G}_{n}(f_{k}))\rightsquigarrow 
\mathcal{N}_{k}(0,\Gamma (f_{i},f_{j})_{1\leq i,j\leq k}),
\end{equation*}

\bigskip \noindent where 
\begin{equation*}
\Gamma (f_{i},f_{j})=\int \left( f_{i}-\mathbb{P}_{Z}(f_{i})\right) \left(
f_{j}-\mathbb{P}_{Z}(f_{j})\right) d\mathbb{P}_{Z}(x),1\leq i ,j\leq k.
\end{equation*}
\end{lemma}

\Bin Now, we recall the results of \cite{Zeghoudi3} on which the previous tool applies.

\subsection{Moments estimators} \label{ssec_02_02}

We have

\begin{proposition} \label{moment} Let $F$ be the Pseudo-Lindley \textit{cdf} corresponding to the \textit{pdf} $f$ in \eqref{glind}. Then

$$
\mathbb{E}(X^k)=\frac{k! (\beta+k)}{\theta^k\beta}, \ k\geq 1,
$$ 

\Bin and so

$$
\mathbb{E}(X)=\frac{\beta+1}{\theta\beta}  \ and \ \mathbb{V}ar(X)=\frac{\beta^2+2\beta-1}{\theta^2 \beta^2}.
$$
\end{proposition}

\Bin To prepare the moment estimation, we recall the empirical means and variances, for $n\geq 1$,

$$
\overline{X}_n=\frac{1}{n}\sum_{j=1}^{n} X_j \ and \ S_n^2=\frac{1}{n}\sum_{j=1}^{n} (X_j-\overline{X}_n)^2.
$$

\Bin Since $F$ has two parameters, the moment estimators of parameters $\theta$ and $\beta$ are solutions of the system of two equations

$$
\mathbb{E}(X)=\overline{X}_n  \ and \ \mathbb{V}ar(X)=S_n^2.
$$

\Bin The solution is given below.

\begin{proposition} \label{estimMoment}  The moments estimators of the parameters $\theta$ and $\beta>1$ for a random sample of size $n$ from the Pseudo-Lindley distribution are given by

$$
\hat{\beta}_n=\frac{\overline{X}_n\sqrt{2}-\sqrt{\overline{X}_n^2-S_n^2}}{\sqrt{\overline{X}_n^2-S^2_n}}
$$

\Bin and

$$
\hat{\theta}_n=\frac{\sqrt{2}}{\overline{X}_n \sqrt{2}-\sqrt{\overline{X}_n^2-S_n^2}}.
$$
\end{proposition}

\Bin \textbf{Remark}. The estimators are the same as those given by \cite{Zeghoudi3}, which are

$$
 \frac{S_n^2+\overline{X}_n^2}{\overline{X}_n^2-S_n^2+\overline{X}_n\sqrt{2} \sqrt{\overline{X}_n^2-S_n^2}} \ (for \ \beta), \  and \ \frac{2\overline{X}_n+\sqrt{2}\sqrt{\overline{X}_n^2-S_n^2}}{\overline{X}_n^2+S_n^2} (for \ \theta).
$$

\Bin  Since \cite{Zeghoudi3} did not show computations and we want to work with our simpler expressions from the point of view of our expansions below, we give the computations that led  to them. $\Diamond$\\

 \Ni \textbf{Proof of Proposition \ref{estimMoment}}. We notice that for all $\beta \in \mathbb{R}$,

$$
\beta^2+2\beta -1 = (\beta+1)^2-2.
$$

\Bin Using that, the moment estimation equations become : $\overline{X}_n=(\beta+1)/(\theta \beta)$  and $S_n^2=((\beta+1)^2-2)/(\theta \beta)^2$.  So we have

$$
(\theta \beta)^2=\frac{(\beta+1)^2}{\overline{X}_n^2}=\frac{(\beta+1)^2-2}{S_n^2}.
$$

\Bin Putting $t^2=(\beta+1)^2$, we get

$$
t=\beta+1=\frac{\overline{X}_n\sqrt{2}}{\sqrt{\overline{X}_n^2-S_n^2}},
$$

\Bin leading to

$$
\beta=\frac{\overline{X}_n\sqrt{2}-\sqrt{\overline{X}_n^2-S_n^2}}{\sqrt{\overline{X}_n^2-S_n^2}}.
$$

\Bin Next, we have 

$$
\theta=\frac{\beta+1}{\overline{X}_n\beta}=\left(\frac{1}{\overline{X}_n} \frac{\overline{X}_n\sqrt{2}}{\sqrt{\overline{X}_n^2-S_n^2}}\right) \left(\frac{\sqrt{\overline{X}_n^2-S_n^2}}{\overline{X}_n\sqrt{2}-\sqrt{\overline{X}_n^2-S_n^2}}\right),
$$

\Bin which gives

$$
\theta=\frac{\sqrt{2}}{\overline{X}_n\sqrt{2}-\sqrt{\overline{X}_n^2-S_n^2}}. \ \blacksquare
$$

\Bin Next, we give an asymptotic joint multivariate normality results.

\section{Asymptotic Joint Multinormal Law of the pair of moments estimators}\label{ssec_03}

\begin{theorem} \label{glind-Law} Suppose that $X$ follows Pseudo-Lindley law of parameters $\beta>1$ and $\theta>0$. Let us denote

$$
m=\mathbb{E}(X)=\frac{\beta+1}{\theta \beta} \ and \ \sigma^2=\mathbb{V}ar(X)=\frac{(\beta+1)^2-2}{(\theta \beta)^2}.
$$ 

\Bin Let $\eta^2=m^2-\sigma^2=2(\theta\beta)^{-2}$ and $\lambda=m\sqrt{2}-\eta$. Define the functions $h_{i}(z)=z^i$ of  $z\in \mathbb{R}$, for $i\geq 1$ and the functions $H_i(z)$ of $z\in \mathbb{R}$, $i \in\{1,2\}$:

$$
H_1=\frac{2}{\eta \lambda}h_1-\frac{1}{\eta \lambda^2\sqrt{2}}h_2;
$$

$$
H_2=-\frac{\lambda(\eta\sqrt{2}+2m)}{\eta^3}h_1+\frac{\lambda+\eta}{2\eta^3}h_2.
$$

\Bin Then, as $n\rightarrow +\infty$,

$$
\biggr(\sqrt{n}(\hat{\theta}_n -\theta), \ \sqrt{n}(\hat \beta_n -\beta)\biggr) \rightsquigarrow \mathcal{N}_2(0, \Sigma),
$$

\Bin where $\rightsquigarrow$ stands for the convergence in distribution, $\mathcal{N}_2(0, \Sigma)$ is a bivariate centered Gaussian random variable with variance-covariance matrix defined by

$$
\Sigma_{11}=\mathbb{V}ar(H_1(X)), \ \Sigma_{22}=\mathbb{V}ar(H_2(X)) \ and \ \Sigma_{12}=\Sigma_{21}=\mathcal{C}ov(H_1(X),H_2(X)).
$$
\end{theorem}

\Bin Before we give the proof, let us see how to compute the covariance matrix. It is proved in \cite{Zeghoudi3} that $X \sim \mathcal{L}_{ps}$ has all moments and

$$
\forall k\geq 1, \ \mathbb{E}X^k=\frac{k!(\beta+k)}{\theta^k\beta}.
$$

\Bin  The $H_i$'s are of the form $H_i=a_ih_1+b_ih_2$, $i \in \{1,2\}$. So we have for $i \in \{1,2\}$,

\begin{equation}
\gamma_i=\mathbb{E}H_i(X)=a_i \frac{(\beta+1)}{\theta\beta} +2 b_i \frac{\beta+2}{\theta^2\beta}. \label{sigma1}
\end{equation}

\Bin Next, let

\begin{equation}
\tau_i^2=\mathbb{E}H_i(X)^2=2a_i^2 \frac{(\beta+2)}{\theta^2\beta} + 24b_i^2 \frac{\beta+4}{\theta^4\beta} + 12 a_i b_i \frac{\beta+3}{\theta^3\beta} \label{sigma2}
\end{equation}

\Bin and

\begin{equation}
c=\mathbb{E}(H_1(X)H_2(X))^2=2 a_1b_1 \frac{\beta+2}{\theta^2\beta}+ 6(a_1b_2+b_1a_2) \frac{\beta+3}{\theta^3\beta}+ 24 a_2b_2 \frac{\beta+4}{\theta^4\beta}. \label{sigma3}
\end{equation}

\Bin We have the following result.

\begin{corollary} \label{corHalimi} Given  $X \sim \mathcal{L}_{ps}$ and the functions $H_i$'s, $i\in \{1,2\}$, in Theorem \ref{glind-Law}, the bi-dimensional random
vector $(H_1(X),H_2(X)$ has the covariance matrix defined by

$$
\Sigma_{ii}=\mathbb{V}ar(H_i(X))=\tau_i^2-\gamma_i^2,\ 1\leq i\leq 2, \ and \ \Sigma_{12}=\mathbb{C}ov(H_1(X),H_2(X))=c-(\gamma_1 \gamma_2),
$$

\Bin where the $\gamma_i$'s, the $\tau_i$'s and $c$ are previously defined in Formulas \eqref{sigma1}, \eqref{sigma2} and \eqref{sigma3}.
\end{corollary}

\Bin \textbf{Proof}. To prove (a), we use the method described in \cite{logs}, itself based on the functional empirical process described above. The method of Subsection 3.1 and Lemma 2 are repeatedly used below. The manipulations of infinitely small quantities in probability (the $o_{\mathbb{P}}$'s) and the bounded quantities in probability (the $O_{\mathbb{P}}$'s) are intensively used. We denote $h_{\ell}(x)=x^{\ell}$, $\ell\geq 1$.\\

\noindent In particular we note that for all $h$ satisfying \eqref{var}, we have 

\begin{equation}
\overline{h(X)}_n=\frac{1}{n}\sum_{1\leq j \leq n} h(X_j)=\mathbb{E}h(X) + \frac{1}{n} \mathbb{G}_n(h),
\end{equation}

\Bin where $\mathbb{G}_n(h)=O_{\mathbb{P}}(1)$ as $n\rightarrow +\infty$. Also for all $h$ satisfying \eqref{var}, for any real-valued function $g : \mathbb{R}\rightarrow \mathbb{R}$ continuous at $\mathbb{E}h(X)$, the application of the delta method to the equation above leads to,

\begin{equation}
g(\overline{h(X)}_n)=g(\mathbb{E}h(X)) + \frac{1}{n} \mathbb{G}_n(g^{\prime}(\mathbb{E}h(X))h) + o_{\mathbb{P}}(n^{-1/2}),
\end{equation}

\Bin as $n\rightarrow +\infty$. We also have (see Lemma 2 in \cite{logfep-tool})

\begin{lemma} \label{lemma1}
Let ($A_{n})$ and ($B_{n})$ be two statistics of the samples $(X_1,\cdots,X_n)$, $n\geq 1$. Let $A$ and $B$ be two real numbers and let $L(z)$ and $H(z)$ be two measurable real-valued functions of  $z\in S$. Suppose that $A_{n}=A+n^{-1/2}\mathbb{G}_{n}(L)+o_{\mathbb{P}}(n^{-1/2})$ and $A_{n}=B+n^{-1/2}\mathbb{G}_{n}(H)+
o_{\mathbb{P}}(n^{-1/2})$. Then

\begin{equation*}
A_{n}+B_{n}=A+B+n^{-1/2}\mathbb{G}_{n}(L+H)+o_{\mathbb{P}}(n^{-1/2}),
\end{equation*}

\Bin
\begin{equation*}
A_{n}B_{n}=AB+n^{-1/2}\mathbb{G}_{n}(BL+AH) + o_{\mathbb{P}}(n^{-1/2})
\end{equation*}

\Bin and if $B\neq 0$,

\begin{equation*}
\frac{A_{n}}{B_{n}}=\frac{A}{B}+n^{-1/2}\mathbb{G}_{n}\left(\frac{1}{B}L-\frac{A}{B^{2}}H\right)+o_{\mathbb{P}}(n^{-1/2}).
\end{equation*}
\end{lemma}

 \Bin Let us now apply this. We denote $m=\mathbb{E}(X)$, $m_2=\mathbb{E}(X^2)$ and $\sigma^2=m_2-m^2$. All the expansions below hold as $n\rightarrow+\infty$ or for a fixed $n\geq 1$. We also use the functions  $h_{\ell}(x)=x^{\ell}$ of $x \in \mathbb{R}$, $\ell\geq 1$. We have
	
%(m^2-\sigma^2)
%\overline{X}_n^2-S_n^2
%\sqrt{\overline{X}_n^2-S_n^2}

%sqrt{m^2-\sigma^2}
%\frac{1}{\sqrt{n}}\mathbb{G}_n\left(\right)

\begin{equation*} \label{Zfor1} 
\overline{X}_n=m + \frac{1}{\sqrt{n}}\mathbb{G}_n(h_1),
\end{equation*}

\Bin which leads, by the delta-method, to

\begin{equation} \label{Zfor2}
\overline{X}_n^2=m^2 + \frac{1}{\sqrt{n}}\mathbb{G}_n(2m h_1) + o_{\mathbb{P}}(n^{-1/2}).
\end{equation}

\Bin Also

\begin{equation*} \label{Zfor3} 
S_n^2=\frac{n}{n-1} \overline{X}_n^2= \frac{n}{n-1} \biggr(m_2+\frac{1}{\sqrt{n}}\mathbb{G}_n(h_2) - \overline{X}_n^2\biggr).
\end{equation*}

\Bin By using the fact that 

$$
C_n=\left(m_2+\frac{1}{\sqrt{n}}\mathbb{G}_n(h_2)\right)-\overline{X}_n^2=O_{\mathbb{P}}(1)
$$ 

\Bin and $n/(n-1)=(1 - 1/(n-1))$, we get

\begin{equation*} \label{Zfor4} 
S_n^2 =  m_2+\frac{1}{\sqrt{n}}\mathbb{G}_n(h_2) - \overline{X}_n^2+ o_{\mathbb{P}}(n^{-1}),
\end{equation*}

\Bin which, by linearity of the functional empirical process and by using $m_2-m^2=\sigma^2$, gives

\begin{equation} \label{Zfor5}
S_n^2= \sigma^2 + \frac{1}{\sqrt{n}}\mathbb{G}_n\left(h_2 -2m h_1 \right) + o_{\mathbb{P}}(n^{-1/2}).
\end{equation}

\Bin We recall that $\eta^2=m^2-\sigma^2$ and $\lambda=m\sqrt{2}-\eta$. Putting together Formulas \eqref{Zfor2} and \eqref{Zfor5} yields

\begin{equation} \label{Zfor6}
\overline{X}_n^2-S_n^2= \eta^2+ \frac{1}{\sqrt{n}}\mathbb{G}_n\left(4mh_1 -h_2) \right) + o_{\mathbb{P}}(n^{-1/2}).
\end{equation}

\Bin Next, by the delta-method, 

\begin{equation} \label{Zfor7}
\sqrt{\overline{X}_n^2-S_n^2}= \eta + \frac{1}{\sqrt{n}} \mathbb{G}_n \left(\frac{4mh_1-h_2}{2\eta} \right) + o_{\mathbb{P}}(n^{-1/2}).
\end{equation}

\Bin By continuing with $\eta^2=m^2-\sigma^2$, we have

\begin{equation} \label{Zfor8} 
\overline{X}_n\sqrt{2} - \sqrt{\overline{X}_n^2-S_n^2}=m\sqrt{2}-\eta + \frac{1}{\sqrt{n}}\mathbb{G}_n\left(\frac{(2\sqrt{2}\eta-4m)h_1-h_2}{2\eta}\right) + o_{\mathbb{P}}(n^{-1/2}).
\end{equation}

\Bin Finally, by the delta-method using $g(z)=1/z$ of $z \in \mathbb{R}\setminus \{0\}$,  we reach the following expansion of $\hat \theta_n$, 

\begin{equation*} \label{Zfor9}
\hat \theta_n= \frac{\sqrt{2}}{\lambda} + \frac{1}{\sqrt{n}}\mathbb{G}_n\left(\frac{\sqrt{2}(4m-2\eta \sqrt{2})}{2\eta \lambda^2}h_1+\frac{\sqrt{2}}{2\eta \lambda^2}h_2\right) + o_{\mathbb{P}}(n^{-1/2}).
\end{equation*}

\Bin Since $\eta^2=m^2-\sigma^2=2/(\theta \beta)^2$, we have

$$
\frac{\sqrt{2}}{m\sqrt{2}-\sqrt{m^2-\sigma^2}}=\frac{\sqrt{2}}{\left(\sqrt{2}\frac{\beta+1}{\theta \beta}-\frac{\sqrt{2}}{\theta \beta}\right)},
$$

\Bin we conclude that

\begin{equation*} 
\hat \theta_n= \theta + \frac{1}{\sqrt{n}}\mathbb{G}_n\left(\frac{4m-2\eta \sqrt{2}}{\eta \lambda^2\sqrt{2}}h_1+\frac{1}{\eta \lambda^2\sqrt{2}}h_2\right) + o_{\mathbb{P}}(n^{-1/2}).
\end{equation*}

\Bin Since $4m-2\eta \sqrt{2}=2\sqrt{2}\lambda$, we finally get

\begin{equation} \label{Zfor9}
\hat \theta_n= \theta + \frac{1}{\sqrt{n}}\mathbb{G}_n\left(\frac{2}{\eta \lambda}h_1-\frac{1}{\eta \lambda^2\sqrt{2}}h_2\right) + o_{\mathbb{P}}(n^{-1/2}).
\end{equation}

\Bin The computations concerning $\hat \theta_n$ are complete. As to that $\hat {\beta}_n$, we use the second expansion in Lemma \ref{lemma1} to Formulas \eqref{Zfor7} and \eqref{Zfor8} with 

$$
A=\lambda, \ L=\frac{2\eta\sqrt{2}-4m}{2\eta}h_1+\frac{1}{2\eta}h_2, \ B=\eta \ and \ H=\frac{4m}{2\eta}h_1-\frac{1}{2\eta}h_2.
$$

\Bin We have

\begin{eqnarray*}
&&\frac{1}{\eta}L - \frac{\lambda}{\eta^2}H\\
&&=\biggr(\frac{2\eta \sqrt{2}}{2\eta^2}h_1+\frac{1}{2\eta^2}h_2\biggr)-\biggr(\frac{4\lambda m}{2\eta^3}h_1-\frac{\lambda}{2\eta^3}h_2\biggr)\\
&&=\frac{\eta(2\eta \sqrt{2}-4m)-4m\lambda}{2\eta^3}h_1+\frac{\lambda+\eta}{2\eta^3}h_2.
\end{eqnarray*}

\Bin We also have $A=\theta \sqrt{2}$ and $B=\eta=\sqrt{2}/(\theta \beta)$ and hence $A/B=\beta$ and so we get

\begin{equation*} \label{Zfor10}
\hat{\beta}_n= \beta + \frac{1}{\sqrt{n}}\mathbb{G}_n\left(\frac{\eta(2\eta \sqrt{2}-4m)-4m\lambda}{2\eta^3}h_1+\frac{\lambda+\eta}{2\eta^3}h_2\right) + o_{\mathbb{P}}(n^{-1/2}).
\end{equation*}

\Bin But $\eta(2\eta \sqrt{2}-4m)-4m\lambda=-2\lambda(\eta \sqrt{2}+2m)$ and finally we have

\begin{equation} \label{Zfor10}
\hat{\beta}_n - \beta = \frac{1}{\sqrt{n}}\mathbb{G}_n\left(-\frac{\lambda(\eta\sqrt{2}+2m)}{\eta^3}h_1+\frac{\lambda+\eta}{2\eta^3}h_2\right) + o_{\mathbb{P}}(n^{-1/2}).
\end{equation}

\Bin Putting together Formula \eqref{Zfor9} and \eqref{Zfor10} leads to

$$
\biggr(\sqrt{n}(\hat{\beta}_n -\beta), \ \sqrt{n}(\hat \theta_n -\theta)\biggr)=(\mathbb{G}_n(H_1), \mathbb{G}_n(H_2))+o_{\mathbb{P}}(1).
$$

\Bin The final conclusion comes from the multivariate asymptotic law of the functional empirical process as given in Lemma \ref{lemma2}. $\blacksquare$

\section{Simulation study} \label{ssec_04}

\subsection{Performance of the moment estimators} \label{ssec_03_01}

For generating data from $X \sim \mathcal{L}_{ps}(\theta, \beta)$, $\theta>0$, $\beta>1$, we use the \textit{cdf}  

$$
1-F(x)= \frac{(\beta+\theta x)e^{-\theta x}}{\beta}, \ x\geq 0,
$$

\Bin and use the quantile function $F^{-1}$ to get samples through 

$$
X_n=F^{-1}(U_n), \ n \geq 1,
$$

\Bin  where $\{U_1, U_2, \cdots \}$ are \textsl{iid} $(0,1)$-uniform random variables. Here, $F^{-1}$ is both the inverse function and the generalized inverse function defined as  

$$
F^{-1}(u)=\inf \{x \in (0,1), \ F(x)\geq u\}.
$$

\Bin  Although $F^{-1}$ does not have an explicit form, we may compute $F^{-1}(u)$, $u \in ]0,1[$, by the dichotomy method. The VB6 codes for the \textit{cdf} and the quantile functions are given in \cite{loSPVB6}. In the section \textit{New and non-classical probability laws, Part 1}, the codes for the \textit{cdf}, the quantile function and data generation are given. In the section \textit{New and non-classical empirical statistical procedures, Part 1}, codes for functions $H_i$ ($i\in \{1,2\}$) and the simulation package both for the convergence and the statistical tests are provided. The results below were performed by those packages.\\

\noindent For fixed $n$, $\theta$ and $\beta$, $B=1000$ iterated generations of $n$ \textit{iid} data from $X \sim \mathcal{L}_{ps}(\theta, \beta)$  and mean values
of the estimators (MVE) $\hat \theta_n$ and $\hat{\beta}_n$ and the square roots of the related Mean Squared Errors (MSE) are computed. The results reported in Table \ref{tab1} for different values of $n$ and for $\theta=\beta=2$.

\begin{table}
	\centering
		\begin{tabular}{lllllll}
		\hline \hline
		Sizes &  $\hat \theta_n$ & $\hat \beta_n$ & MSE($\hat \theta_n$) & MSE($\hat \beta_n$) & p-value ($\theta$)& p-value ($\beta$)\\
		\hline \hline
		n=50	& 2.17	& 2.18	& 0.45& 1.9 & 6.5\%& 4.5\%\\
		200		&2.04 	& 2.17	& 0.23& 1.61 & 5.7\%& 4.05\%\\
		375		&2			& 2.13	& 0.18& 0.988 & 4.1\% & 3.85\%\\
		400		&2.01		& 2.16	& 0.17& 0.87 & 3.81\% & 3.61\%\\
		500		&2 			& 2.06	& 0.16& 0.9 & 3.20\%& 3.61\%\\
		700		&2 			& 2.07	& 0.13& 0.47 & &\\
		900		&2			& 2.03	& 0.11& 0.4 & &\\
		1000	&2			& 2.08	& 0.1& 0.36 & &\\
		1500	&2			& 2.01	& 0.09& 0.26 & &\\
		2000	&2			& 2.02	& 0.08& 0.22 & &\\
		2500	&2			& 2.01	& 0.07& 0.2 & &\\
		\hline \hline
		\end{tabular}
	\caption{convergence of the mean values estimators and the square roots of the MSE's versus sample sizes. P-values of statistical tests in the last two columns}
	\label{tab1}
\end{table}

\begin{figure}[h]
	\centering
		\includegraphics[width=0.90\textwidth]{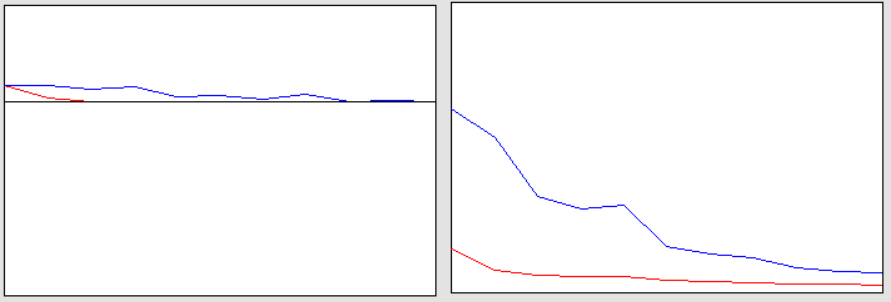}
	\caption{Left: convergence of the mean values versus sample sizes. Right: convergence of the square roots of the MSE's versus sample sizes}
	\label{fig01}
\end{figure}

\Bin  \textbf{Comments}. In the current illustration $\theta=2$ and $\beta=2$, we have a very good performance only for sizes greater than $400$. For sizes of that order and greater, the mean value error and the mean squared error are outstanding for the estimation of $\theta$, far better than for the estimation of $\beta$ which can be considered as good enough. In Figure \ref{fig01}, we see the different performance in terms of mean squared errors (MSE) and the not so different values of mean value errors. How do we explain that the performance is lower for $\beta$'s estimation? Indeed, we have 

$$
\sqrt{m^2-\sigma^2})=1/(\theta\beta).
$$

\Bin  The current value is $\sqrt{m^2-\sigma^2}=0.5$. $\hat{\beta}_n$ has $\sqrt{\overline{X}_n^2-S_n^2}$ as denominator. So a weaker value of $\sqrt{m^2-\sigma^2}$ will induce a weaker value of $\sqrt{\overline{X}_n^2-S_n^2}$ that will have effected to make $\beta$ bigger. But if $\sqrt{\overline{X}_n^2-S_n^2}$ is very weak, the denominator of $\hat{\theta}_n$, which is $\overline{X}_n\sqrt{2}-\sqrt{\overline{X}_n^2-S_n^2}$ is more stable around $\overline{X}_n\sqrt{2}$ and this entails that $\hat{\theta}_n$ is stable for a significant expectation $m$.\\

\noindent The recommendation for application is the following. Given the size of the data, we compute the estimation of $\theta$ and $\beta$ and next proceed to the simulations for these parameters. If the results are good enough, we make the hypothesis $X \sim \mathcal{L}_{ps}(\theta,\beta)$ and apply the statistical test below.\\

\subsection{Statistical tests} \label{sec_03_01} $ $

\Ni The columns  \textit{p-value ($\theta$)} and \textit{p-value ($\beta$)} in Table \ref{tab1} are enough to show that the statistical tests based on the asymptotic laws are accurate around and above sample sizes of $n=50$.\\

\section{Conclusion and perspectives}\label{ssec_05} Our contribution is that the provided statistical tests are accurate for medium sample sizes and can be used in Reliability and in Survival Analysis. But in medical sciences, especially when blood is concerned, the data may be scarce. Investigations aiming at obtained accurate tests for small sizes and based on some nonparametric methods are under way. 

%\frac{1}{\sqrt{n}}\mathbb{G}_n\left(\right)
%\begin{equation*} \label{Zfor1} 
%\end{equation*}
%\begin{equation*} \label{Zfor1} 
%\end{equation*}
%\begin{equation*} \label{Zfor1} 
%\end{equation*}
%\begin{equation*} \label{Zfor1} 
%\end{equation*}
%\begin{equation*} \label{Zfor1} 
%\end{equation*}
%\begin{equation*} \label{Zfor1} 
%\end{equation*}
%\begin{equation*} \label{Zfor1} 
%\end{equation*}
%\begin{equation*} \label{Zfor1} 
%\end{equation*}

%\label{exampleRmetric}
%\begin{definition} \label{topo_02_} 
%\end{definition}
%\begin{proposition} \label{topo_02_} 
%\end{proposition}
%\begin{lemma} \label{topo_02_} 
%\end{lemma}
%\begin{eqnarray} 
%\end{eqnarray}
%\begin{equation} 
%\end{equation}
%\begin{corollary} \label{topo_02_} 
%\end{corollary}
%\begin{theorem} \label{topo_02_} 
%\end{theorem}

\include{lindleyGeneralizedByHalim_appendix}
\end{document}